\documentclass[12pt]{article}
\usepackage{amsmath,amssymb,amsthm,amsfonts,mathtools,mathrsfs}
\usepackage{etaremune, enumerate, float, verbatim}
\usepackage{algorithm}
\usepackage{algorithmic}
\usepackage{hyperref, theoremref}
\usepackage{soul} 
\usepackage{graphicx}
\usepackage{url}
\usepackage[mathlines]{lineno}
\usepackage{dsfont} 
\usepackage{tikz, graphicx, subcaption, caption}
\usepackage[algo2e,ruled,vlined]{algorithm2e}
\usepackage{xcolor}
\numberwithin{figure}{section}
\usepackage{theoremref}
\usepackage{listings}
\usepackage[noadjust]{cite}
\usetikzlibrary{positioning}

\newtheorem{thm}{Theorem}[section]
\newtheorem{cor}[thm]{Corollary}
\newtheorem{lem}[thm]{Lemma}
\newtheorem{prop}[thm]{Proposition}

\newtheorem{obs}[thm]{Observation}

\theoremstyle{definition}

\theoremstyle{definition}

\theoremstyle{definition}
\newtheorem{ex}[thm]{Example}

\newcommand{\R}{\mathbb{R}}

\newcommand{\mptn}{\mathcal{S}} 

\DeclareMathOperator{\rank}{rank}

\DeclareMathOperator{\vspan}{span}
\DeclareMathOperator{\nul}{nul}
\DeclareMathOperator{\col}{col}
\DeclareMathOperator{\spark}{spark}
\DeclareMathOperator{\mr}{mr}
\newcommand{\mrplus}{\mr_{+}}
\DeclareMathOperator{\supp}{supp}
\DeclareMathOperator{\GN}{GN}
\DeclareMathOperator{\GM}{GM}
\DeclareMathOperator{\mfsr}{mfsr}
\newcommand{\mfsrplus}{\mfsr_{+}}
\DeclareMathOperator{\Block}{Block}
\DeclareMathOperator{\failed}{F}

\title{
The Spark of Symmetric Matrices Described by a Graph}
\date{\today}
\author{Louis Deaett\thanks{Department of Mathematics and Statistics, Quinnipiac University, Hamden, CT, USA. louis.deaett@quinnipiac.edu}
\and Shaun Fallat
\thanks{Department of Mathematics and Statistics, University of Regina, Regina, Saskatchewan, Canada. shaun.fallat@uregina.ca}
\and Veronika Furst \thanks{Department of Mathematics, Fort Lewis College, Durango, CO, USA. furst\_v@fortlewis.edu}
\and
John Hutchens \thanks{Department of Mathematics and Statistics, University of San Francisco, San Francisco, CA, USA. jhutchens@usfca.edu}~\thanks{Corresponding Author} \and Lon Mitchell \thanks{Department of Mathematics \& Statistics, Eastern Michigan University, Ypsilanti, MI, USA.  lmitch50@emich.edu} 
\and
Yaqi Zhang\thanks{Department of Mathematics, Drexel University, Philadelphia, PA, USA. yaqizhangus@outlook.com}
}

\begin{document}
\maketitle\vspace{-20pt}

\begin{abstract} 
We investigate the sparsity of null vectors of real symmetric matrices whose off-diagonal pattern of zero and nonzero entries is described by the adjacencies of a graph.  We use the definition of the spark of a matrix, the smallest number of nonzero coordinates of any null vector, to define the spark of a graph as the smallest possible spark of a corresponding matrix.
We study connections of graph spark to well-known concepts including minimum rank, forts, orthogonal representations, Parter and Fiedler vertices, and vertex connectivity.
\end{abstract}

\noindent {\bf Keywords:} Null vectors, maximum nullity, spark, zero forcing, forts, connectivity, generic nullity, minimum rank.
\medskip

\noindent {\bf AMS subject classification:} 05C50, 15A18 (primary) 15A29 (secondary).

\medskip


\section{Introduction}
Denote the set of all real symmetric $n \times n$ matrices by $S_n(\mathbb{R})$, and suppose $A=[a_{ij}] \in S_n(\mathbb{R})$. 
We say $G(A)$ is the graph of $A$  if $G(A)$ has the vertex set $V=\{v_1,v_2,\dotsc,v_n\}$ and edge set $E=\{v_iv_j \mid a_{ij} \neq 0, \ i \neq j\}$. Note that $G(A)$ is independent of the values of the diagonal entries of $A$. On the other hand, if $G$ is a graph of order $n$ (i.e., $|G| = |V(G)| =n$) with vertex set \(\{v_1,v_2,\dotsc,v_n\}\), then the set of real symmetric matrices described by the graph $G$ is given by \(\mathcal{S}(G) = \{A \in S_n(\mathbb{R}) \mid G(A) = G \}\).  Here and in what follows, we consider only simple, undirected graphs \(G = (V(G), E(G))\). One of the most captivating and unresolved problems associated with the class $\mathcal{S}(G)$ is the so-called inverse eigenvalue problem for graphs, abbreviated as IEP-G (see \cite{AIMMINRANK, MR4284782, MR3752611, MR2033476, loop-zero-forcing, MR3291662, HLS2022inverse}). This fundamental problem asks for a complete description of the possible spectra realized by the set $\mathcal{S}(G)$ for a given graph $G$. The IEP-G has garnered significant attention over the past 30 years with many fascinating advances and applications (see, for example, the books \cite{HLS2022inverse, MR3752611} and the references therein). However, a complete general resolution is still very much elusive. Notwithstanding this, researchers have developed a wealth of results, implications, and applications tied to the IEP-G (see \cite{IEPG2} for a recent example). In particular, a number of related concepts and parameters have been explored and have shed light on different aspects of the IEP-G\@.  The \emph{minimum rank} of a graph $G$ is defined to be \(\mr(G) = \min \{\rank(A) \mid A \in \mptn(G)\}\). The \emph{maximum nullity} (or \emph{maximum corank}) of a graph $G$ is defined to be $M(G) = \max \{\nul(A) \mid A \in \mathcal{S}(G)\} = n - \mr(G)$, where $\nul(A)$ denotes the nullity of $A$ or the dimension of the null space of $A$, written as $N(A)$.  The {\em minimum semidefinite rank} $\mrplus(G)$ is defined analogously as the minimum rank over all positive semidefinite matrices in $\mptn(G)$.  (We refer the reader to the works \cite{Barrett2004, MR3665573, Nylen, MR1902112, MR3313038, MR2594483}.)  The column space of $A$ will be denoted $\col(A)$.

While the primary focus on the IEP-G has been on the potential list of eigenvalues of matrices in $\mptn(G)$, there is also justified interest in studying the associated eigenvectors or zero/nonzero patterns of the associated eigenvectors. One of the earliest results along these lines is by Fiedler \cite{F2} where the eigenvectors of matrices associated with connected acyclic graphs (or trees) were studied. One by-product of this work was the realization that investigating the zero coordinates of an eigenvector leads to certain implications about a graph (or in the case of \cite{F2} a tree). Since Fiedler's pioneering work in 1975, research into the possible patterns of eigenvectors for matrices associated with a graph has been developed along a number of lines, including: nodal domains, Laplacian eigenvectors (e.g, Fiedler vectors), and more recently zero forcing on graphs (see also \cite{Nylen, nylen1, FN, MR3034535}).  We note here that it is sufficient to study the zero/nonzero patterns of null vectors of $A \in \mptn(G)$, since any eigenvector $x$ corresponding to the eigenvalue $\lambda$ of $A$ can be considered as a null vector of the matrix $A - \lambda I \in \mptn(G)$.

As our work relies heavily on the theory of graphs, we list some useful notation and provide some relevant terminology here before we discuss zero forcing and spark for graphs.
A \emph{subgraph} $H = (V(H), E(H))$ of $G = (V(G),E(G))$ is a graph with $V(H) \subseteq V(G)$ and $E(H) \subseteq E(G)$; $H$ is an {\em induced subgraph} of $G$ if $E(H) = \{vw \in E(G) \mid v, w \in V(H)\}$. 
The {\em complement} of $G = (V,E)$ is the graph $\overline{G} = (V, \overline{E})$, where $\overline{E}$ consists of all pairs of vertices in $V$ that are not contained in $E$.   
We say two vertices $v,w$ are {\em adjacent}, or are {\em neighbors}, if $vw \in E$, and we may write this as $v \sim w$.  Let $N_G(v) = \{w\in V \mid w \sim v\}$ be the open neighborhood of $v$ and denote its cardinality by $\deg(v) = |N_G(v)|$. The closed neighborhood of $v$ is $N_G[v] = N_G(v) \cup \{v\}$.  The minimum degree of the graph is $\delta= \delta(G)= \min\{\deg(v) \mid v \in V(G) \}$.

A \emph{path} is a graph, denoted $P_n$, with $V = \{v_1,\dotsc,v_n\}$, where $v_1,\dotsc,v_n$ are distinct, and $E=\{v_iv_{i+1} \mid i = 1,\dotsc, n - 1\}$. 
A \emph{cycle} $C_n$ on $n$ vertices 
    $V = \{v_1,\dotsc,v_n\}$ has $E = \{v_iv_{i+1} \mid i = 1,\dotsc, n - 1\} \cup \{v_nv_1\}$.
A graph is \emph{connected} if for every pair of distinct vertices $v$ and $u$ there is a path from $v$  to $u$ (and thus also from $u$ to $v$).
A \emph{tree} is a connected graph with no cycles.  A {\em complete graph} $K_n$ on $n$ vertices has $E = \{v_iv_j \mid i\neq j\}$.  A {\em complete bipartite graph} $K_{m,n}$ has vertex set $V = V_1 \cup V_2$, where $|V_1| = m$ and $|V_2|=n$, and edge set \(E = \{v_iv_j \mid v_i\in V_1, \ v_j\in V_2\}\).  If $G=(V(G), E(G))$ and  $H=(V(H), E(H))$ are two graphs, then the {\em Cartesian product} of $G$ and $H$, denoted by $G \mathbin{\Box} H$, is the graph with vertex set $V(G) \times V(H)$ and two vertices $(u,v)$ and $(w,z)$ are adjacent in $G \mathbin{\Box} H$ if and only if $u=w$ and $vz \in E(H)$ or $uw \in E(G)$ and $v=z$.

{\em Zero forcing} is a coloring process involving the vertices of a graph. At the beginning of the process, each vertex is either blue or white, and each type of zero forcing follows a specific color change rule which can change the color of a white vertex to blue. The process stops when no more vertices can be colored blue.  The \emph{standard zero forcing} color change rule is to change the color of a white vertex $w$ to blue if $w$ is the unique white neighbor of a blue vertex $v$. If an initial subset of blue vertices can, after repeated application of the color change rule, change all vertices to blue, then this subset is referred to as a {\em zero forcing set}. Zero forcing was introduced to provide a combinatorial upper bound for $M(G)$ and, in particular, detects subsets of coordinates of a null vector $x$ of any $A\in \mptn(G)$ that, if designated as zero, imply $x$ must in fact be the zero vector. As such, it seems natural to study the zero coordinates in null vectors (see \cite{AIMMINRANK, loop-zero-forcing, z-paper1, HLS2022inverse} for more details). 

More precisely, given a real vector $x$, the \emph{support} of $x$ is the collection of indices $i$ for which $x_i \neq 0$. We denote the support of $x$ by $\supp(x)$. Suppose $A \in \mptn(G)$ and $Ax=0$. A basic consequence of the zero forcing process outlined above is that if $\supp(x)$ is disjoint from a zero forcing set for $G$, then $x=0$. 

Suppose $B \subseteq V$ is initially colored blue, and that $B'$ is the set of all blue vertices obtained from $B$ by repeatedly applying the color change rule. We call $B'$ the closure of $B$.  If nonempty, the subset $V \setminus B'$ (the remaining white vertices) is known as a fort in $G$. In fact, a \emph{fort} in a graph is a nonempty subset $F$ of vertices such that no vertex outside of $F$ is adjacent to exactly one vertex of $F$ (see \cite{brim}). 

Forts are naturally connected to the support of null vectors. As we are interested in sparse null vectors, we seek to determine forts of minimum size in a given graph. Finally, it is a simple observation in basic linear algebra that if $x$ is in $N(A)$, for any matrix $A$, then the columns of $A$ that correspond to $\supp(x)$ must form a linearly dependent set. This leads us to the notion of the spark of a matrix, which we present in the next section. 

This paper is organized into sections combining various topics with the spark of a graph. In Section 2, we define the spark of a graph and explore a connection with forts in the graph. In Section 3, we discuss relationships between the concepts of spark and rank. Then in Section 4, we investigate an association between spark and the vertex connectivity of a graph, and we generalize a theorem concerning orthogonal representations of graphs. In Section 5, we pay particular attention to graphs with small spark, and we close our work with some further connections in Section 6.

\section{Spark and forts of graphs} \label{spark&forts}

As our main focus is studying the support of null vectors, and, in particular, to exhibit null vectors that have small support, we begin with the notion of the spark of a matrix. Namely, the \emph{spark} of a matrix $A$ is the smallest integer $s$ such that there exists a set of $s$ columns in $A$ which are linearly dependent, i.e., $\spark(A)$ is the minimum size of the support of a null vector of $A$.  If $A\in S_n(\R)$ is nonsingular, $\spark(A)$ is defined to be $n+1$.  Sparse solutions to underdetermined linear systems, and thereby the concept of spark, have gained significant attention in compressed sensing (see \cite{MR1963681, MR2241189, Compressed-sensing}).  Computing the spark of a matrix is known to be NP-hard \cite[Problem A6.MP5]{MR519066}.
We define the \emph{spark} of a graph $G$ to be 
\[ \spark(G)=\min_{A\in \mathcal{S}(G)} \spark(A) . \]
Note that, for every graph $G$, the Laplacian matrix of $G$ gives a singular matrix in $\mptn(G)$, showing that $\spark(G) \le n$.
In addition, it is not hard to see that $\spark(G)=1$ if and only if $G$ contains an isolated vertex.  Furthermore, if $G$ is disconnected, then $\spark(G)$ is obtained by simply minimizing the spark across all of the connected components of $G$. Thus, we assume from this point on that all graphs considered are connected and hence contain no isolated vertices.

We illustrate the above notions with the following example.

\begin{ex}
Let $G$ be a graph on 5 vertices consisting of a 5-cycle on vertices $\{1,2,3,4,5\}$ with two additional edges $13$ and $25$. Suppose $A \in \mathcal{S}(G)$ is given by 
\[A = \left[\begin{matrix}
1 &1 & 1& 0 & 1\\ 
1 &1 & 1& 0 & 1\\
1 &1 & 3& 1 & 0\\
0 &0 & 1& 3 & 1\\
1 &1 & 0& 1 & 3\\
\end{matrix}\right] \; {\rm and} \;\
x=\left[\begin{matrix}
1\\ -1\\0 \\ 0 \\ 0 \end{matrix}\right].\] 
Observe that $Ax=0$, and hence $\spark(A)\leq 2$. Since $G$ is connected (or more precisely has no isolated vertices) it is clear that $\spark(G)>1$. Thus $\spark(A) = \spark(G) = 2$. Finally we note that the pair $\{1,2\}$ forms a fort in $G$.
\end{ex}

The connection between the support of  a null vector $x$ for some $A \in \mptn(G)$
and a fort in $G$ in the previous example is a known result (although it may not be published); we provide a proof here for completeness, as our primary aim is studying the support of null vectors.  Recall that the columns of $A\in \mptn(G)$ are indexed by the vertices of $G$, and we use the column indices and the graph vertices interchangeably. 

\begin{thm}\label{FortSupport}
    For any matrix $A \in \mptn(G)$, the support of any nonzero null vector of $A$ is a fort of $G$. Conversely, for any fort $F$ of $G$ and any vector $x$ whose support is $F$, there is a matrix $A \in \mptn(G)$ that has $x$ as a null vector.  That is, $\spark(G)$ is the cardinality of a minimum fort of $G$.
\end{thm}
\begin{proof}
Given a matrix $A\in \mptn(G)$ and a vector $x\neq 0$ such that $Ax=0$, let $W = \{j \in V(G) \mid j \in \supp(x)\}$.  Suppose there exists $i\in V(G) \setminus W$ with exactly one neighbor $j$ in $W$.  Then 
\[ 0 = [Ax]_i = a_{ij} x_j,\]
where $a_{ij}\neq 0$.  Thus, $x_j = 0$, contradicting $j\in \supp(x)$.  So $W$ is a fort of $G$.

Conversely, assume $F$ is a fort of $G$ such that $F = \{i\in V(G) \mid i\in\supp(x)\}$ for some nonzero vector $x$.  We construct $A$ by performing the following steps: 
    \begin{enumerate}
    \item First let $A$ be the adjacency matrix of $G$.  In the next two steps, we modify certain nonzero entries of $A$.
    \item Let $S=\{i \mid x_i=0\} = V(G)\setminus F$. For $i\in S$,  let $B_i= N_G(i) \setminus S = N_G(i) \cap F$. Note that $|B_i|\neq 1$ by the definition of a fort. For $j\in B_i$ and $j\neq \max B_i$, set $A[i,j]=A[j,i]=1/x_j$; if $j=\max B_i$, then set $A[i,j]=A[j,i]=(1-|B_i|)/x_j$.  
    \item For $k\in \supp(x)$,  assign $A[k,k]=-\frac{[Ax]_k}{x_k} = -\frac{\sum_{j\neq k} a_{kj}x_j}{x_k}$.
    \end{enumerate}   
    For $i\in F$, 
    \[
    [Ax]_i=a_{ii}x_i+\sum_{j\neq i} a_{ij}x_j=-\frac{\sum_{j\neq i} a_{ij}x_j}{x_i}x_i+\sum_{j\neq i} a_{ij}x_j=0. 
    \]
    For $i\in V(G)\setminus F$, 
    \begin{align*}
    [Ax]_i &= \sum_{j\nsim i}a_{ij}x_j+a_{ii}x_i+\sum_{j\sim i}a_{ij}x_j=\sum_{j\sim i}a_{ij}x_j\\
    &= \sum_{j\in N_G(i)\cap S}a_{ij}x_j+\sum_{j\in N_G(i)\setminus S} a_{ij}x_j = \sum_{j\in B_i}a_{ij}x_j\\
    &= \sum_{\substack{j\in B_i \\ j\neq \max{B_i}}} \frac{1}{x_j}x_j+\frac{1-|B_i|}{x_{\max B_i}}x_{\max B_i}=0.\qedhere
    \end{align*}
\end{proof}

Although $\spark(G)$ is defined in reference to the matrices in $\mathcal{S}(G)$, Theorem \ref{FortSupport} shows that in fact this parameter can be defined entirely in graph-theoretic terms.  That is, the spark of a graph does not have to be defined in terms of the spark of any matrices.

The next two propositions explore possible sizes of forts of graphs in more detail.

\begin{prop}\label{newallforts}
Let $G$ be a graph with minimum degree $\delta$.  Then every subset $W \subseteq V(G)$ with $|W| = n-m +1$ is a fort of $G$ if and only if $m \leq \delta$.
\end{prop}

\begin{proof}
Assume $m\leq \delta$, and consider $W \subseteq V(G)$ with $|W| = n-m +1 \geq n-\delta + 1$.  For any vertex $v \notin W$, there are at most $\delta -2$ vertices not in $W \cup \{v\}$.

Conversely, suppose $m > \delta$ and $v_0 \in V(G)$ such that 
$ N_G(v_0) =\{ v_1,v_2, \dotsc, v_{\delta}\}$;
     we can then label \[ V(G) =\{ v_0,v_1, \dotsc, v_{\delta}, s_{\delta +2}, \dotsc, s_n \}. \]
    Since $m+1\geq \delta + 2$, the set 
    $ W = \{v_1, s_{m+1}, \dotsc, s_n \} $
    satisfies $|W|=n-m+1$ but is not a fort.
\end{proof}

\begin{prop} \label{everyfortlarger}
If every $k$-subset of $V(G)$ is a fort of $G$, then every $(k+1)$-subset of $V(G)$ is a fort of $G$.
\end{prop}
\begin{proof}
Let every $k$-subset of $V(G)$ be a fort of $G$, and let $W \subseteq V(G)$ such that $|W|=k+1$.  In particular $|W|\geq 2$.  Assume $W$ is not a fort, so there exists $x \in V(G)\setminus W$ such that $N_G(x) \cap W =\{w_1\}$, where
\[ W = \{w_1,w_2, \dotsc, w_{k+1} \}. \]
Then $W' = W\setminus \{w_2\}$ is not a fort since $x \in V(G) \setminus W'$ and $N_G(x) \cap W' = \{w_1\}$.  But $|W'|=k$, giving us a contradiction.
\end{proof}

In light of Proposition \ref{everyfortlarger}, it is interesting to note that simply adding vertices to a fort does not guarantee that the resulting set is a fort.  We present examples of graphs that skip fort sizes after Theorem \ref{spark_anyrank}.

Related to the notion of a fort is the notion of a \emph{failed zero forcing set} \cite[Definition 1.4]{Fetcie_failed_forcing}.  This is simply a subset of vertices that is not a zero forcing set; however, it is interesting to ask for the largest size of such a set.  This is known as the \emph{failed zero forcing number} $\failed(G)$ of the graph $G$~\cite{Fetcie_failed_forcing}.  
The complement of a failed zero forcing set has also been called a \emph{zero blocking set}, with the smallest size of such a set called the \emph{zero blocking number}, and denoted by $\Block(G)$ \cite{zeroblocking}.
As noted in \cite{zeroblocking}, a set is a failed zero forcing set of maximum size if and only if its complement is a fort of minimum size. Using Theorem \ref{FortSupport}, it follows that $\Block(G)$ and $\spark(G)$ are the same.

In summary, we have the following observation.

\begin{obs}
Let $G$ be a graph on $n$ vertices.  Then $\Block(G)=\spark(G)$ and $\failed(G) = n-\spark(G)$.
\end{obs}

Hence, the problem of determining the failed zero forcing number of a graph and the problem of determining its zero blocking number are both equivalent to determining the spark of the graph.
In fact, this problem, like that of computing the spark of a matrix, is NP-hard \cite{Shitov2017}.

\section{Spark and rank of matrices associated with a graph}

The spark and rank of a matrix \(A \in \mptn(G)\) are clearly related, as the definitions give directly that \(\spark(A) \leq \rank(A) + 1\). In this section, we investigate when this inequality becomes an equality.  We say a matrix has \emph{full spark} if \(\spark(A) = \rank(A) + 1\).  Analogously to $\mr(G)$ and $\mrplus(G)$, we define the \emph{minimum full spark rank} of a graph $G$ as $\mfsr(G) = \min \{ \rank(A) \mid A \in \mptn(G), \ \rank(A) = \spark(A) -1 \}$ and $\mfsrplus(G)$ as the corresponding minimum full spark rank for positive semidefinite matrices.  

The next result is a core result in linear algebra and lays the groundwork for establishing a relationship between spark and rank.

\begin{thm}\label{PSMnssparkkplus1}
Let \(A\) be a symmetric \(n \times n\) real matrix with \(\rank(A) = k \).  If \(k=n\) then each \(k \times k\) principal submatrix of \(A\) is nonsingular and \(A\) has full spark.  

If \(k < n\) and \(X\) is an \(n \times (n-k)\) real matrix with \(\rank(X) = n-k\) such that \(AX=0\), then the following are equivalent:
\begin{enumerate}[(1)]
\item Each \(k \times k\) principal submatrix of \(A\) is nonsingular.
\item Each \((n-k) \times (n-k)\) submatrix of \(X\) is nonsingular.
\item \(\spark(A)=k+1\).
\end{enumerate}
\end{thm}

\begin{proof} If \(k=n\) the result is clear, so assume \(k < n\).

$(1) \Leftrightarrow (2)$: Without loss of generality we can write \(A = \left[\begin{matrix}B & C^T \\ C & M \end{matrix}\right]\) and \(X=\left[\begin{matrix} Y \\ Z\end{matrix}\right] \), where \(B\) is \(k\times k\) and \(Z\) is  \((n-k) \times (n-k)\), and argue that \(B\) is singular if and only if \(Z\) is singular.  

If \(Z\) is singular, then there exists a nonzero vector \(v\) with \(Zv = 0\).  Since \(\rank(X)=n-k\), \(Yv\) must be nonzero, and so  \(Xv\) is a nonzero vector in the nullspace of \(A\).  Then \(Yv\) is a nonzero vector in the nullspace of \(B\), so \(B\) is singular.

If \(B\) is singular, then there exists a nonzero vector \(v\) with \(v^TB=0\).  If \(v^TC^T \neq 0\), then \(v^TC^TZ= 0\) implies we are done.  So assume \(v^TC^T = 0\).  Then \(\left[\begin{matrix} B \\ C\end{matrix}\right]v=0 \), so that \(\left[\begin{matrix} v \\ 0\end{matrix}\right] \) is a nonzero vector in the nullspace of \(A\).  Since \(\rank(X)=n-k\), the columns of \(X\) are a basis for the nullspace of \(A\).  Thus the vector \(\left[\begin{matrix} v \\ 0\end{matrix}\right] \) is a nontrivial linear combination of the columns of \(X\), and so \(Z\) is singular.

\((1)\Rightarrow (3)\): If each \(k \times k\) principal submatrix of \(A\) is nonsingular, then each set of \(k\) columns of $A$ is linearly independent and \(\spark(A) \geq k+1\); that is, \(A\) has full spark.

\((3)\Rightarrow (1)\): Assume \(\spark(A) = k+1\), and suppose that there exists a \(k\times k\) principal submatrix $B$ of $A$ that is singular.  Without loss of generality, write \(A = \left[\begin{matrix}B & C^T \\ C & M \end{matrix}\right]\) in block form.  Let \(X=\left[\begin{matrix} Y \\ Z\end{matrix}\right] \) in similar block form be a matrix whose columns form a basis for the nullspace of \(A\), so that \(AX=0\).  By the proof of $(1) \Leftrightarrow (2)$, the \((n-k) \times (n-k)\) matrix \(Z\) must also be singular, and there exists a nonzero vector \(w\) with \(Zw= 0\).  But then \(  \left[\begin{matrix} Y \\ Z\end{matrix}\right]w = \left[\begin{matrix} Yw \\ 0\end{matrix}\right] \) is a nonzero vector (if \(Yw=0\), then the columns of \(X\) are linearly dependent) in the nullspace of \(A\), so that \(BYw=CYw = 0\). But then \(\left[\begin{matrix} B \\ C\end{matrix}\right] Yw = 0 \)  implies that the columns of \(\left[\begin{matrix} B\\ C\end{matrix}\right]\) are linearly dependent, which contradicts \(\spark(A)=k+1\).   
\end{proof}

We next consider a bordering-type result concerning the spark of a symmetric matrix.
\begin{lem}
Suppose $A$ is an $n \times n$ real symmetric matrix. Consider the bordered $(n+1) \times (n+1)$ symmetric matrix given by
\(B = \left[\begin{matrix}x^T Ax & x^T A \\ Ax & A \end{matrix}\right]\) for some vector $x$.
Then
\begin{enumerate}[(1)]
    \item \(\rank(B)=\rank(A)\);
    \item if \(\lvert\supp(x)\rvert =k\), then \(\spark(B) \leq k+1\).
    \end{enumerate}
\end{lem}
\begin{proof}
Statement (1) is trivial. For (2), observe that the vector 
\(\left[\begin{matrix} -1 \\ x\end{matrix}\right]\) is a null vector for $B$, and the result follows.
\end{proof}
A simple consequence of the above lemma can be deduced if we assume in addition that $A$ is invertible. Then $\rank(B)=n$, and thus it follows that \(\spark(B)=\lvert\supp(x)\rvert +1\), since the dimension of the null space of $B$ is one.

Given a graph $G$ with order $n$, we are interested in finding all possible ordered pairs of integers $(k, s)$, $1\leq k\leq n$ and $2\leq s\leq n+1$, such that there exists $A\in\mptn(G)$ with $\rank(A) = k$ and $\spark(A) = s$.  Note that $k=n$ if and only if $s = n+1$.  After finding a matrix with some fixed spark $s$ and minimum corresponding rank $k\geq s-1$, the next result shows that all higher ranks are achievable with the same spark. Before we state this result, we recall the following observation:  for any symmetric matrix $A$, the $j$th standard basis vector $e_j \in \col(A)$ if and only if $j \not\in \supp(x)$ for all $x \in N(A)$, which follows easily from the fact that the null space of a symmetric matrix $A$ is the orthogonal complement of the column space of $A$.

\begin{thm} \label{spark_anyrank}
If $A \in \mptn(G)$ such that $\rank(A)=k < n-1$ and $\spark(A)=s$, then there exists a matrix 
$B \in \mptn(G)$ such that $\rank(B)=k+1$ and $\spark(B) =s$.
\end{thm}
\begin{proof}
Assume $A \in \mptn(G)$ such that $\rank(A)=k < n-1$ and $\spark(A)=s$.  Let us choose a basis for the null space of $A$
\[ N(A) = \vspan\{ \eta_1,\eta_2, \eta_3, \dotsc, \eta_{n-k} \}, \]
such that $\lvert\supp(\eta_1)\rvert =s$.  

Let $\eta_i = (y_{i1},y_{i2}, \dotsc, y_{in})^T$ for $1 \leq i \leq n-k$.  We claim there is a matrix $B = A +D \in \mptn(G)$ where $D$ is a diagonal matrix, such that $\eta_1 \in N(B)$, $\eta_2 \not\in N(B)$ and $N(B) \subset N(A)$.  Since $\lvert\supp(\eta_1)\rvert =s$ we know $\supp(\eta_1) \neq \supp(\eta_2)$, otherwise there would be a null vector whose support size is smaller than $s$.  

So we can choose $j \in \supp(\eta_2) \setminus \supp(\eta_1)$.  Choose a new basis for $N(A)$ such that
\[ N(A) = \vspan\{ \eta_1,\eta_2, \eta_3', \dotsc, \eta_{n-k}' \}, \]
with $\eta_i' = \frac{y_{ij}}{y_{2j}} \eta_2 - \eta_i$ for $3 \leq i \leq n-k$.   Then  $j \not \in \supp(\eta_i')$ for $3 \leq i \leq n-k$ and $j \not \in \supp(\eta_1)$.   

Let $e_j$ represent the $j$th standard basis vector in $\R^n$ and consider $e_j e_j^T = E_{jj}$.  We see that 
\[ (A + E_{jj})\eta = A \eta + E_{jj} \eta =E_{jj} \eta = e_j e_j^T \eta \]
for $\eta \in N(A)$.  Notice that $e_j \not \in \col(A)$ since $e_j^T \eta_2 = y_{2j} \neq 0$ and $\col(A) = N(A)^{\perp}$.  By \cite{MR311674}, this implies $\rank(A+E_{jj}) = \rank(A) +\rank(E_{jj}) = k+1$.  

Now $\eta_2 \not\in N(A+E_{jj})$ but 
\[\{ \eta_1, \eta_3', \dotsc, \eta_{n-k}'\} \subset N(A+E_{jj}) \]
since each vector in the set is orthogonal to $e_j$.  This gives us
\[ N(A+E_{jj}) = \vspan\{ \eta_1, \eta_3', \dotsc, \eta_{n-k}' \} \]
since this is a set of $n-k-1$ linearly independent vectors in $N(A+E_{jj})$ where $\dim(N(A+E_{jj})) = n-k-1$.  Since $N(A+E_{jj}) \subset N(A)$, $\lvert\supp(\eta_1)\rvert =s$, and $\eta_1 \in N(A+E_{jj})$, we have that $\spark(A+E_{jj}) = s$.
\end{proof}

We note here that given $A \in \mptn(G)$ we cannot necessarily find another matrix $B \in \mptn(G)$ such that $\spark(B)=\spark(A)+1$ and $\rank(B)=\rank(A)$.  This follows, in part, due to the fact that if $G$ has a fort of size $s$ this may not guarantee that $G$ has a fort of size $s+1$ or $s-1$.  Define the \emph{fort sequence} of $G$ to be the sequence of the form $(s_2,s_3,\dotsc, s_n)$ where $G$ has $n$ vertices and $s_i$ is the number of forts in $G$ with $i$ vertices.  Note that we are beginning the fort sequence at $s_2$; since we only consider graphs without isolated vertices, all graphs we consider have $s_1=0$. There are many examples of graphs that skip fort sizes.  

For example, consider a \emph{spider graph} (also known as a {\em generalized star}), which is a tree with one vertex having degree greater than $2$, the \emph{central vertex}, and all other vertices having degree at most $2$.  The paths radiating out from the central vertex are called the \emph{legs} and do not contain the central vertex.  We denote such a graph as $sp(n_1,n_2,\dotsc,n_l)$ where $l$ is the degree of the largest-degree vertex (i.e., the number of legs) and $n_j$ is the number of vertices in each leg.  So the order of $sp(n_1,n_2,\dotsc,n_l)$ is $1+\sum n_j$.  

\begin{figure}[ht!]
\begin{center}
\begin{tikzpicture}[
fillednode/.style={circle, draw=black!100, fill=gray!75, thick, scale=0.75},
unfillednode/.style={circle, draw=black!100, fill=black!0, thick, scale=0.75},
]
\node[unfillednode,label={$a$}]      (a)                            {}  ;
\node[unfillednode,label={$c$}]        (c)       [below right= of a] {};
\node[unfillednode,label={$b$}]        (b)       [below left= of c] {};
\node[unfillednode,label={$v_1$}]        (v1)       [right=1 cm of c] {};
\node[unfillednode,label={$v_2$}]        (v2)       [right=1 cm of v1] {};
\node[unfillednode,label={$v_{m-1}$}]        (vm-1)       [right=1 cm of v2] {};
\node[unfillednode,label={$v_m$}]        (vm)       [right=1 cm of vm-1] {};

\draw[thick] (a) -- (c);
\draw[thick] (b) -- (c);
\draw[thick] (c) -- (v1);
\draw[thick] (v1) -- (v2);
\draw[thick,dashed] (v2) -- (vm-1);
\draw[thick] (vm-1) -- (vm);

\end{tikzpicture}
\end{center}
\caption{Spider graph $sp(m,1,1)$}\label{spm11}
\end{figure}
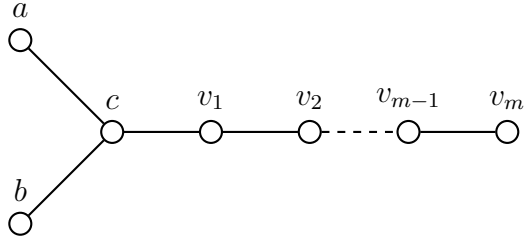

Now consider a special class of spider graphs of the form $sp(m,1,1)$, where $m>3$, depicted in Figure \ref{spm11}.  

The smallest fort size is $2$ corresponding to the unique minimum fort $\{a,b\}$.  Any fort $F\subseteq V$ with $|F| \geq 3$ must contain the vertex $v_m$, otherwise $v_{i+1}$ where $i = \max\{j<m \mid v_j\in F\}$ (here considering $c$ as $v_0$) would be adjacent to only one vertex in $F$.
The next smallest fort is a minimum fort for $P_{m+2}$, which arises as the induced subgraph of $G$ on either $\{a,c,v_1,\ldots,v_m\}$ or $\{b,c,v_1,\ldots,v_m\}$.  The path $P_{m+2}$ has a minimum fort of size $\lceil\frac{m+3}{2} \rceil$; note that the minimum fort will not contain $c$ for either parity of $m$.
So the fort sequence for $sp(m,1,1)$ is of the form $(1,0, \ldots, 0, s_{\lceil\frac{m+3}{2} \rceil}, \ldots, s_{m+3})$ with $s_i\neq 0$ for $\lceil\frac{m+3}{2} \rceil \leq i \leq m+3$.

Moreover, if $sp(m,1,\dotsc, 1)$ has $3\leq l<m$ legs, then the fort sequence is of the form $(s_2,\ldots,s_{l-1},0 \ldots, 0, s_{\lceil\frac{m+3}{2} \rceil}, \ldots, s_{m+l})$, where $s_i\neq 0$ for $2 \leq i \leq l-1$ or $\lceil\frac{m+3}{2} \rceil \leq i \leq m+l$.  Hence there is no bound on size of the gap between two nonzero fort sizes in a graph or constraints on where in the sequence a gap can occur.

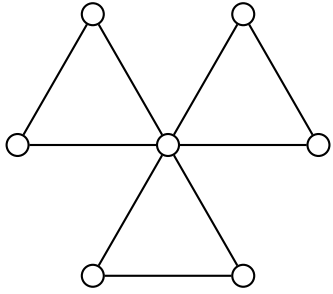
\begin{figure}[ht!]
\begin{center}
\begin{tikzpicture}[
fillednode/.style={circle, draw=black!100, fill=gray!75, thick, scale=0.75},
unfillednode/.style={circle, draw=black!100, fill=black!0, thick, scale=0.75},scale=2
]
\node[unfillednode] (v0) at (0,0)  {};
\node[unfillednode] (v1) at (0.5,0.87)  {};
\node[unfillednode] (v2) at  (1,0)  {};
\node[unfillednode] (v3) at (0.5,-0.87)  {};
\node[unfillednode] (v4) at (-0.5,-0.87) {};
\node[unfillednode] (v5) at (-1,0)  {};
\node[unfillednode] (v6) at  (-0.5,0.87)  {};
\draw[thick] (v0) to (v1);
\draw[thick] (v0) to (v2);
\draw[thick] (v0) to (v3);
\draw[thick] (v0) to (v4);
\draw[thick] (v0) to (v5);
\draw[thick] (v0) to (v6);
\draw[thick] (v1) to (v2);
\draw[thick] (v3) to (v4);
\draw[thick] (v5) to (v6);
\end{tikzpicture}
\end{center}
\caption{Friendship graph $F_3$}\label{F_3}
\end{figure}

An example of a graph that is not a tree and skips a fort size is the friendship graph $F_3$ shown in Figure \ref{F_3}. This graph has forts of size $2$ and $4$ but no forts of size $3$.  Indeed, any pair of adjacent non-central vertices forms a fort, but any set $S$ of three vertices in $V=V(F_3)$ must leave at least one vertex in $V\setminus S$ adjacent to only one vertex in $S$; any pair of adjacent pairs of non-central vertices forms a fort of size $4$.  In fact $F_3$ has fort sequence $(3,0,11,12,7,1)$.

\section{Spark and connectivity of graphs}

The {\em vertex connectivity} of a graph $G$, denoted by $\kappa(G)$, is defined as the minimum size of a set of vertices whose deletion disconnects the graph. Such a set of vertices is known as a {\em cut set}. Further, we say a graph is {\em $k$-connected} if $\kappa(G) \geq k$. 

For a graph \(G\), a \emph{(faithful) orthogonal representation} of \(G\) of dimension \(k\) is a set of vectors in \(\mathbb{R}^k\), one corresponding to each vertex, with the property that two vertices are nonadjacent if and only if their corresponding vectors are orthogonal. An orthogonal representation of \(G\) in \(\mathbb{R}^k\) is in {\em general position} if every subset of \(k\) vectors is linearly independent.  Note that this is equivalent to the existence of a positive semidefinite matrix $A\in\mptn(G)$ with $k = \rank(A) = \spark(A) - 1$; \(A\) is called the \emph{Gram matrix} of the orthogonal representation~\cite{horn_johnson_2012}.

\begin{thm}[\cite{connectivity,correction}]\label{thm:LSS}
For a graph \(G\) with \(n\) vertices, the following are equivalent:
\begin{enumerate}[(1)]
    \item \(G\) is \((n - k)\)-connected.
    \item \(G\) has a general position orthogonal representation in \(\mathbb{R}^k\).
    \item \(G\) has an orthogonal representation in \(\mathbb{R}^k\) consisting of unit vectors such that for each vertex \(v\) the vectors representing the vertices not adjacent to \(v\) are linearly independent.
\end{enumerate} 
\end{thm}

A consequence of Theorem \ref{thm:LSS} is that the minimum semidefinite full spark rank is dictated by the connectivity of the graph, with \(\mfsrplus(G) = n-\kappa(G)\) for every graph \(G\).  Indeed, if the rank of a positive semidefinite matrix drops below this threshold, then the spark may be forced to drop even further (recall that, in general, \(\delta(G) \geq \kappa(G)\)):

\begin{cor}\label{cor:fourtwo}
If $A$ is a positive semidefinite matrix for a graph $G$ with vertex connectivity $\kappa(G)$ and $\rank A < n - \kappa(G)$, then $\spark A \leq n - \delta(G) -1$.
\end{cor}

\begin{proof}
If $\spark A > n - \delta(G) - 1$, then every set of $n - \delta(G) - 1$ vertices is linearly independent.  Since every vertex has at most that many non-neighbors, every set of non-neighbors is linearly independent, which implies $G$ is $(n-\rank(A))$-connected by Theorem \ref{thm:LSS}.
\end{proof}

The minimum semidefinite full-spark rank of a graph may be strictly larger than the minimum semidefinite rank, as demonstrated by the following example.  Let $G = C_4 \mathbin{\Box} P_t$ be the Cartesian product of the cycle $C_4$ and the path $P_t$ on $t\geq 2$ vertices. The minimum semidefinite rank of $G$ is $4t-4$~\cite{psdMandZ}, but   $\delta(G) = 3 = \kappa(G)$, so any positive semidefinite matrix in \(\mptn(G)\) with rank \(4t-4\) cannot have full spark, and the smallest possible rank of a full spark positive semidefinite matrix in $\mptn(G)$ is $4t-3$.  That is, \(4t-3 = \mfsrplus(G) > \mrplus(G) = 4t-4\).  

The minimum rank and minimum semidefinite rank of $G = C_4 \mathbin{\Box} P_t$ coincide, with $\mr(G) = 4t-4$~\cite{AIMMINRANK}, so we can ask if there exists a full spark symmetric matrix for $G$ that has minimum rank but is not positive semidefinite.  Perhaps surprisingly, we show in our next result that it is not possible to achieve a lower full-spark rank with arbitrary symmetric matrices. That is, \(\mfsr(G) = \mfsrplus(G) = n-\kappa(G)\) for every graph \(G\).

\begin{thm}\label{prop:extension}
A graph \(G\) is \((n-k)\)-connected if and only if there exists \(A \in \mptn(G)\) with \(k = \rank(A) = \spark(A)-1\).  
\end{thm}

\begin{proof}
One direction follows from Theorem \ref{thm:LSS}.  For the other direction, let \(A \in \mptn(G)\) with \(k = \rank(A) = \spark(A)-1\) and suppose that \(G\) is not \((n-k)\)-connected.  Then there exists a cut set \(\alpha\) of \(n-k-1\) vertices that leaves at least two connected components and we can write the matrix \(A(\alpha)\) where the rows and columns corresponding to \(\alpha\) are removed in block form as
\[
A(\alpha) = \left[  \begin{matrix} B & 0 \\ 0 & C \end{matrix}\right].
\]
Since \(A\) has rank \(k\), \(A(\alpha)\) has rank at most \(k\) but size \((k+1) \times (k+1)\), so that \(A(\alpha)\) must be singular.  Without loss of generality, that means that \(B\) must also be singular.  Since \(B\) is at most a \(k \times k\) matrix, \(A(\alpha)\), and thus \(A\), contains a singular \(k \times k\)  principal submatrix (any \(k \times k\) submatrix of \(A(\alpha)\) that has \(B\) as a submatrix will also be singular because of the block structure).  But this contradicts Theorem \ref{PSMnssparkkplus1}.   
\end{proof}

In light of Theorem \ref{prop:extension}, it is natural to ask if all of Theorem \ref{thm:LSS} could extend to arbitrary symmetric matrices.  While Theorem \ref{PSMnssparkkplus1} tells us \(A \in \mptn(G)\) with \(k=\rank(A)= \spark(A) - 1\) has each \(k \times k\) submatrix invertible, smaller submatrices need not be invertible.  For example, 
    \[
    \begin{bmatrix}
     0 & 0 & \phantom{-}3 & 1 & 4\\
     0 & 2 & \phantom{-}4 & 4 & 4\\
     3 & 4 & -4 & 0 & 0\\
     1 & 4 & \phantom{-}0 & 4 & 0\\
     4 & 4 &  \phantom{-}0 & 0 & 8
    \end{bmatrix}
    \]
is a rank-three matrix in \(\mptn(K_{2,3})\) (for the bipartite graph \(K_{2,3}\), note that \(\kappa(K_{2,3}) = 2\))  with every $3\times 3$ principal submatrix nonsingular but with singular principal submatrices of sizes two and one.  In particular, the $1\times 1$ principal submatrix corresponding to vertex 1, the non-neighbor of one of the degree-three vertices, is singular. Thus we cannot extend Theorem \ref{thm:LSS} by focusing on principal submatrices; however, we do find a full generalization by looking instead (in the spirit of spark) at linearly independent columns. 

\begin{thm}\label{thm:generalized}
For a graph \(G\) with \(n\) vertices, the following are equivalent:
\begin{enumerate}[(1)]
    \item \(G\) is \((n - k)\)-connected.
    \item There exists \(A \in \mptn(G)\) with \(k = \rank(A) = \spark(A)-1\).
    \item There exists \(A \in \mptn(G)\) with \(k = \rank(A)\) and such that for any vertex \(v\) of \(G\) the columns of \(A\) corresponding to \(v\) and its non-neighbors are linearly independent.
    \item There exists \(A \in \mptn(G)\) with \(k = \rank(A)\) and such that for any vertex \(v\) of \(G\) the columns of \(A\) corresponding to the non-neighbors of \(v\) are linearly independent.
\end{enumerate} 
\end{thm}

\begin{proof}
The equivalence of (1) and (2) is the content of Theorem \ref{prop:extension}, and we will use it to show (2) implies (3).  If (2) is true, then the matrix \(A\) also satisfies (3): by (1) and using \(n-k \leq \kappa(G) \leq \delta(G)\), each vertex \(v\) has at most 
    \[
    n - \delta(G) - 1 \leq n - \kappa(G) - 1 \leq n - (n-k) - 1 = k - 1
    \]
non-neighbors; since \(\spark(A) = k+1\), the at-most-\(k\) columns of \(A\) corresponding to \(v\) and its non-neighbors must be linearly independent.  

Since (3) is stronger than (4), the main work will now be to prove that (4) implies (1).  Suppose \(A \in \mptn(G)\) with \(\rank(A) = k\) is such that for any vertex \(v\) the columns corresponding to the non-neighbors of \(v\) are linearly independent but \(G\) is not \(n-k\) connected.  Then we can find a cut set \(C\) with \(n-k-1\) vertices and can write 
    \[
    A = \left[ \begin{matrix} 
    M_1 & 0 & N_1^T \\
    0& M_2 & N_2^T \\
    N_1 & N_2 & M_3 
    \end{matrix}\right]
    \]
where \(M_3\) corresponds to the vertices of \(C\).  Let each \(M_i\) have size \(d_i \times d_i\) (so \(d_3 = n-k-1\)) and rank \(r_i\).  Finally, let \(n_i =d_i - r_i\) for each \(i\in\{1,2\}\).  We wish to show that \(\rank(A) \geq d_1 + d_2\) in order to get a contradiction.  If \(d_1 = r_1\) and \(d_2 = r_2\), we are done.  So assume without loss of generality that \(n_1 > 0\) and \(n_1 \geq n_2\). By our assumption of (4), the column rank of 
    \[ 
    \left[\begin{matrix}M_1 \\ 0 \\ N_1\end{matrix} \right]
    \] 
is \(d_1\) and the column rank of 
    \[ 
    \left[\begin{matrix}0 \\ M_2 \\ N_2\end{matrix} \right]
    \] 
is \(d_2\).   However, unlike the positive semidefinite case, all we can say is that the column rank of
    \[ 
    \left[\begin{matrix}M_1 & 0 \\ 0 & M_2 \\ N_1 & N_2 \end{matrix} \right]
    \]
is at least \(d_1 + r_2\).  And yet, by symmetry, the row rank of 
    \[ 
    \left[\begin{matrix}M_1 &  0 & N_1^T\end{matrix} \right]
    \] 
is \(d_1\), and thus so is its column rank.  Because $n_1 > 0$, $M_1$ is singular, so there exists a set of $r_1$ columns of $M_1$ and $n_1$ columns of $N_1^T$ that is linearly independent.  

Let $S$ be the set consisting of the first $d_1$ columns of $A$, and let $T$ be the $n_1$ columns among the last $d_3$ columns of $A$ corresponding to the selected columns of $N_1^T$.  Since $S$ is a linearly independent set and the selected $n_1$ columns of $N_1^T$ are not in $\col(M_1)$, $S\cup T$ is a linearly independent set.  Find a basis of $r_2$ vectors for $\col(M_2)$, and let $U$ denote the corresponding columns of $A$.  A linear dependence relation among the vectors in $S\cup U\cup T$ would imply a linear dependence relation among the vectors in $U\cup T$ since the entries in rows $d_1+1$ through $d_1+d_2$ are zeros in each vector in $S$.  Moreover, the entries in rows $1$ through $d_1$ are zeros in each vector of $U$, implying a linear dependence relation among the vectors in $T$.  Each of the sets $S$, $U$, and $T$ is linearly independent, so working backwards, all three linear dependence relations must be trivial, and $S\cup U\cup T$ is linearly independent.  Thus $A$ has at least \(d_1 + r_2 + n_1 \geq d_1 + r_2 + n_2 = d_1 + d_2\) linearly independent columns.

\end{proof}

\section{Graphs with Small Spark}

Recall from Section \ref{spark&forts} that $\spark(G) \geq 2$ for any graph $G$ with no isolated vertices and that $\spark(G)$ is the size of the smallest possible fort in $G$.  The following lemma characterizes graphs $G$ with $\spark(G) = 2$. 

\begin{lem} \label{duplicates}
Let $G$ be a graph. Then $\spark(G)=2$ if and only if there exists $u,v\in V(G)$ such that (1) $uv\in E(G)$ and $N_G[u]=N_G[v]$ or (2) $uv\notin E(G)$ and $N_G(u)=N_G(v)$.
\end{lem}
\begin{proof}

Assume $\spark(G) = 2$.  By Theorem \ref{FortSupport}, $G$ has a minimum fort of size 2, say $F = \{u,v\}$.  By the definition of a fort, every vertex in $V(G)\setminus F$ is adjacent to neither or both of the vertices in $F$, implying either condition (1) or (2).  Conversely, if condition (1) or (2) hold, then $F = \{u,v\}$ is a fort in $G$; having size $2$, $F$ must be a minimum fort.
\end{proof}

If $u,v \in V(G)$ satisfy either condition (1) or condition (2) of Lemma \ref{duplicates}, then we will refer to them as {\em duplicate} vertices.  

\begin{lem} \label{NoDuplicates}
    Let $G$ be a graph of order $n\geq 3$.  Then $G$ must have duplicate vertices if either of the following conditions hold:
    \begin{enumerate}[(1)]
        \item $\mr(G) \leq 2$,
        \item $\kappa(G) \geq n-2$.
    \end{enumerate}
\end{lem}

\begin{proof}
    Note that by Theorem \ref{thm:LSS}, $\mrplus(G) \leq n-\kappa(G)$, so condition (2) implies condition (1).  Therefore, it suffices to prove that (1) implies the existence of duplicate vertices.

 Assume $\mr(G) \leq 2$.  By Theorem 9 of \cite{Barrett2004}, $\overline{G}$ can be expressed as the union of at most two complete graphs and of bipartite graphs. Suppose $\overline{G}$ consists of only complete graphs, in which case it must consist of the union of exactly two complete graphs since $G$ is connected.  Then $G$ is a complete bipartite graph of order $n\geq 3$ and therefore has two duplicate vertices.  On the other hand, suppose $\overline{G}$ has a complete bipartite graph as a component. This component must contain at least three vertices, and we can let $u$ and $v$ be vertices in the same partite set; then $N_{\overline{G}}(u) = N_{\overline{G}}(v)$, so $N_G[u] = N_G[v]$, and $u$ and $v$ are duplicate vertices in $G$.
\end{proof}

\begin{prop}
    If $G$ is a graph with $\spark(G) \geq 3$, then $\mr(G) \geq 3$ and $\kappa(G) \leq n-3$.  In particular, if $\spark(G)=3$ and $A\in \mptn(G)$ with $
    \spark(A) = 3$, then $A$ is not full spark. 
\end{prop}

\begin{proof}
    Since $\spark(G) \geq 3$, $G$ has no no duplicate vertices by Lemma \ref{duplicates}.  The result then follows by Lemma \ref{NoDuplicates}.
\end{proof}

\begin{prop}
If $\spark(G) = 2$, then either $G$ is a path on three vertices or $\mr(G) < n-1$.
\end{prop}

\begin{proof}
If $\spark(G) =2$, then $G$ has a pair of duplicate vertices by Lemma \ref{duplicates}.  If $\mr(G) = n-1$ then $G$ is a path $P_n$ on $n$ vertices \cite{MR2350678}, which can contain duplicate vertices only if $n=3$.   
\end{proof}

Naturally, we can ask if an analogous result holds for graphs with larger spark.  Unfortunately, increasing the spark by 1 does not necessarily decrease the minimum rank's bound by 1, as is illustrated by the following example. 

\begin{ex}
    Let $G = C_n$ be a cycle on $n$ vertices and let $H$ be obtained from $G$ by adding the edges $v_1v_4$ and $v_1v_{n-2}$. Then, for $n \geq 5$, it follows that $H$ has no duplicate vertices, so $\spark(H)  \geq 3$. On the other hand, $\{v_1, v_3, v_{n-1}\}$ forms a fort in $H$. Hence $\spark(H)=3$. Finally, it is not difficult to deduce that $\mr(H)=n-2$.
\end{ex}

We saw in Theorem \ref{prop:extension} that considering matrices in $\mptn(G)$ does not provide an advantage over positive semidefinite matrices in achieving minimum rank and full spark.  We may also consider matrices of minimum rank and minimum spark.  This is not necessarily achievable with a positive semidefinite matrix, as the next example demonstrates.  The hypercube graph \(Q_3 = C_4 \mathbin{\Box} P_2\) has minimum rank $4$. The matrix \(H_3\) given in \cite[p.\ 1636]{AIMMINRANK} for the graph \(Q_3\) has rank \(4\) and spark \(3\).  Since $Q_3$ has no duplicate vertices, $\spark(G) > 2$ by Lemma \ref{duplicates}, so $H_3$ achieves minimum rank and minimum spark for $Q_3$.

\begin{prop}\label{nonsingularfullspark}
If \(A \in \mptn(Q_3)\) is positive semidefinite and \(\rank(A)=4\), then \(\spark(A)=4\).  
\end{prop}

\begin{proof}
We have \(\kappa(Q_3) = \delta(Q_3)= 3\). Thus  \(\spark(A) \leq 4\) by Corollary \ref{cor:fourtwo}.  Since \(Q_3\) has no duplicate vertices, \(\spark(A) > 2\) by Lemma \ref{duplicates}.   Suppose that \(\spark(A)=3\).  Then, in the orthogonal representation corresponding to \(A\), we can find three vectors that are linearly dependent.  That is, the dimension of their span must be at most two.  Consider the subgraph corresponding to the three vectors.  It cannot be complete as \(K_3\) is not a subgraph of \(Q_3\).  If it has no edges, then all three vectors are orthogonal and cannot be linearly dependent.  If there is just one edge, then two of the vectors must be orthogonal to the third, which would make them linearly dependent (in a one-dimensional subspace), contradicting $\spark(A)>2$.  If there are two edges, then two of the vectors must be orthogonal to each other, say \(\vec{v}_1\) and \(\vec{v}_2\), and the third must then be a linear combination of both: \(\alpha\vec{v}_1 + \beta\vec{v}_2\) with \(\alpha\beta \neq 0\).  But then the vector representing the third neighbor (in \(Q_3\)) of the vertex represented by \(\alpha\vec{v}_1 + \beta\vec{v}_2\) would have nonzero dot product with at least one of \(\vec{v}_1\) and \(\vec{v}_2\), a contradiction.
\end{proof}

\section{Further Connections}

Let $A\in\mptn(G)$ and $v\in V(G)$. Denote by $A(v)$ the principal submatrix of $A$ obtained by deleting $v$.  The vertex $v\in V(G)$ is a \emph{Parter vertex} (P-vertex) for $A\in\mptn(G)$ if $\nul(A(v)) = \nul(A) + 1$ ($\Leftrightarrow \rank(A) = \rank(A(v)) + 2)$
(see the original works \cite{P, W} and \cite{KIM2008, JS} for more recent related work on these topics).  A {\em Fiedler vertex} (F-vertex) $v\in V(G)$ for $A\in \mptn(G)$ is a vertex that satisfies  $\nul(A(v)) \geq \nul(A)$.
Both Parter and Fiedler vertices are interconnected with zero coordinates in null vectors: 

\begin{lem}[{\cite[Theorem 2.1]{JS}}]\label{ZerosInNullVec}
Let $A\in\mptn(G)$ and $v\in V(G)$. Then $\nul(A(v)) \geq \nul(A)$ if and only if every null vector of $A$ has a $0$ in the $v$-th coordinate.
\end{lem}

According to Kim and Shader~\cite{KIM2008}, if we partition a singular symmetric matrix \(A\) as \(A=\left[ \begin{matrix} a & x^T \\ x & B \end{matrix}\right]\), then vertex 1 is an F-vertex if and only if  \(\left[\begin{matrix} a \\ x  \end{matrix}\right]\) is not in the column span of  \(\left[\begin{matrix}  x^T \\   B\end{matrix}\right]\); and vertex 1 is a P-vertex if and only if \(x\) is not in the column span of \(B\).  Since a positive semidefinite matrix automatically has the row/column inclusion property \cite{FJ}, a positive semidefinite matrix cannot have a P-vertex.  And if \(A\) is positive semidefinite, then \(a \geq y^TBy\) where \(x=By\), so vertex 1 is a F-vertex if and only if \(a > y^TBy \).  In that case, we can decrease the rank of \(A\) by exactly one if we replace \(a\) with \(y^TBy\).  Thus a positive semidefinite matrix in $\mptn(G)$ of minimum (semidefinite) rank cannot have an F-vertex. 

Instead of just considering the support of a particular null vector, there is also interest in considering the support of the null space. That is, for a matrix $A$ in $S_n(\mathbb{R})$, we define the support of the null space of $A$ as 
\[ \supp N(A) = \{ i \mid x_i\neq 0 \mbox{ for some } x \in N(A)\}.\]

An important well-known fact for matrices in $\mptn(T)$, where $T$ is a tree, is the following:

\begin{prop}[\cite{nylen1, F2}] \label{1-dimN} Suppose $T$ is a tree and $A \in \mptn(T)$.  If $\supp N(A)=V(T)$, then $\dim N(A)=1$.
\end{prop}

The minimum semidefinite rank of a tree $T$ on $n$ vertices is $n-1$.  Hence, the converse of Proposition \ref{1-dimN} states that a matrix $A$ realizing this minimum must have full null support.  The following theorem shows that this in fact holds not just for trees but for all graphs.

\begin{thm}
If a positive semidefinite matrix $A \in \mptn(G)$ has $\rank(A) = \mrplus(G)$, then $\supp N(A) = V(G)$.
\end{thm}

\begin{proof}
A positive semidefinite matrix of minimum (semidefinite) rank cannot have an F-vertex, so no vertex has a zero component in every null vector.
\end{proof}

For trees, full null space support turns out to be equivalent to full spark for singular matrices. 

\begin{thm}\label{prop:whatever}
Let \(T\) be a tree and \(A \in \mptn(T)\) be singular.  Then the following statements are equivalent:
\begin{enumerate}[(1)]
    \item $A$ has full spark: \(\spark(A) = \rank(A)+1\).
    \item $A$ has full null space support, that is,   \(\supp N(A)=V(T)\).
    \item \(A\) does not have a Parter vertex.
\end{enumerate}
\end{thm}

\begin{proof}
$(1) \Leftrightarrow (2)$:  Suppose $k = \rank(A) = \spark(A) - 1$.  By Theorem \ref{prop:extension}, $T$ must be $(n-k)$-connected, which implies $k = n-1$ since $T$ is a tree and $A$ is singular.  Since $\dim N(A) = 1$ and $A$ has full spark, Theorem \ref{PSMnssparkkplus1} implies $\supp N(A) = V(T)$.

Conversely, suppose $\supp N(A) = V(T)$.  By Proposition \ref{1-dimN}, $\dim N(A) = 1$.  This implies that every nontrivial null vector of $A$ has only nonzero entries.  By Theorem \ref{PSMnssparkkplus1}, $A$ must have full spark.

$(2) \Leftrightarrow (3)$ follows from Lemma \ref{ZerosInNullVec}.
\end{proof}

We note here for completeness that if $\mr(G) < \mrplus(G)$, then a minimum rank matrix need not have an F-vertex. Suppose $G=K_{2,3}$. Then $\mr(G)=2$ and $\mr_+(G)=3$. The adjacency matrix of $G$ is a minimum rank (indefinite) matrix that does not have an F-vertex.

In \cite{genericnullity} Hogben and Shader define a real matrix \(X\) to be \emph{generic} if every square submatrix of $X$ is nonsingular.  Then the \emph{generic nullity} of a nonzero $A \in \mathbb{R}^{n\times n}$ is 
\[ \GN(A) = \max\{ k  \mid  X \in \mathbb{R}^{n\times k}, \ AX=0, \ X \text{ is generic} \}, \]
and the \emph{maximum generic nullity} of a graph is 
\[ \GM(G) = \max\{ \GN(A) \mid A \in \mptn(G) \}. \]

For any graph $G$, note that $\GM(G)\geq 1$ since the all-ones vector belongs to the null space of the graph's Laplacian matrix.  We end with an interesting relation between rank, spark, and maximum generic nullity of a graph.

\begin{thm} \label{maxgennul}
If there exists \(A \in \mptn(G)\) such that \(\rank(A) = k\) and \(\spark(A) = k+1\), then \(\GM(G) \geq n-k\).  
\end{thm}

\begin{proof}
By Theorem \ref{PSMnssparkkplus1}, each $k\times k$ principal submatrix of \(A\) is nonsingular.  If \(k=n\), the result is clear, so assume \(k < n\). Let \(X\) be a \(n \times (n-k)\) matrix whose columns form a basis for the nullspace of \(A\), so that \(AX = 0\).  By Theorem \ref{PSMnssparkkplus1}, \(X\) is generic.  Thus \(\GN(A) = n-k\) and \(\GM(G) \geq n-k\). 
\end{proof}

Theorem \ref{prop:extension} also has an implication for generic nullity.  An immediate consequence, by Theorem \ref{maxgennul}, is that $\GM(G) \geq \kappa(G)$ for any graph $G$ (Corollary 4.2 of \cite{genericnullity}).  If the inequality is strict, we can say more:

\begin{cor}
If \(\GM(G) > \kappa(G)\), then any matrix \(A \in \mptn(G)\) with \(\GN(A) = \GM(G)\) satisfies \(\GN(A) < \nul(A)\).  
\end{cor}

\begin{proof}
Suppose $\GM(G) = k$ where $k > \kappa(G)$, and let $A\in\mptn(G)$ with $\GN(A) = k$.  Then there exists a generic matrix $X\in\R^{n\times k}$ such that $AX = 0$, implying $\nul(A) \geq k$.  Suppose $\nul(A) = k$.  Then $\rank(A) = n-k$ and $A$ is full spark by Theorem \ref{PSMnssparkkplus1}.  But this contradicts $\mfsr(G) = n-\kappa(G)$.
\end{proof}

\section*{Acknowledgments}
Shaun M. Fallat was supported in part by an NSERC Discovery Research Grant, Application No.: RGPIN--2019--03934. The research of Yaqi Zhang was partially supported by Simons Foundation grant 355645 and NSF grant DMS 2000037.

This project began as part of the ``Inverse Eigenvalue Problems for Graphs and Zero Forcing'' Research Community sponsored by the American Institute of Mathematics (AIM). We thank AIM for their support, and we thank the organizers and participants for contributing to this stimulating research experience.


\bibliographystyle{plainurl}
\bibliography{refs}

\end{document}